\documentclass[12pt]{amsart}

\usepackage[T1]{fontenc}
\usepackage{lmodern}
\usepackage{amsmath,amssymb,mathtools}
\usepackage{enumitem}
\usepackage{microtype}
\usepackage[colorlinks=true,linkcolor=blue,citecolor=blue,urlcolor=blue]{hyperref}
\hypersetup{
  pdftitle={Non-Group Property, Beltrami Coefficients, and Contractibility of Harmonic Quasiconformal Boundary Maps},
  pdfauthor={Tailiang Liu and Yuliang Shen},
  pdfsubject={Beltrami coefficients, boundary nonclosure, and BMO contractibility of harmonic quasiconformal mappings},
  pdfkeywords={harmonic quasiconformal mapping, Beltrami coefficient, Hilbert transform, Carleson measure, BMO Teichmuller space}
}

\newtheorem{theorem}{Theorem}[section]
\newtheorem{proposition}[theorem]{Proposition}
\newtheorem{lemma}[theorem]{Lemma}
\newtheorem{corollary}[theorem]{Corollary}
\theoremstyle{remark}
\newtheorem{remark}[theorem]{Remark}
\numberwithin{equation}{section}

\newcommand{\D}{\mathbb D}
\newcommand{\Hup}{\mathbb H}
\newcommand{\T}{\mathbb T}
\newcommand{\R}{\mathbb R}

\newcommand{\BMO}{\mathrm{BMO}}
\newcommand{\SQS}{\mathrm{SQS}}
\newcommand{\HQS}{\mathrm{HQS}}
\newcommand{\HQC}{\mathrm{HQC}}
\newcommand{\QS}{\mathrm{QS}}
\newcommand{\HR}{\mathcal H_{\R}}
\newcommand{\HT}{\mathcal H_{\T}}
\newcommand{\id}{\mathrm{id}}
\newcommand{\conv}{\operatorname{conv}}

\title[Harmonic Quasiconformal Maps]
{Harmonic Quasiconformal  Mappings in the BMO-Teichm\"uller space}

\author{Tailiang Liu}
\address{School of Mathematics and Physics, Jiangsu University of Technology, Changzhou 213001, P. R. China}
\email{ltlmath@jsut.edu.cn}

\author{Yuliang Shen}
\address{Department of Mathematics, Soochow University, Suzhou 215006, P. R. China}
\email{ylshen@suda.edu.cn}
\thanks{Project supported by the National Natural Science Foundation of China (Nos.~12401095 and 12571083).}

\subjclass[2020]{30C62, 30F60, 30H35, 42B20}
\keywords{Harmonic quasiconformal mapping, Beltrami coefficient, boundary correspondence, Hilbert transform, BMO-Teichm\"uller space}

\begin{document}

\begin{abstract}
	Let $\HQS(\R)$ denote the class of increasing homeomorphisms of the real line $\mathbb{R}$
	that admit harmonic quasiconformal extensions to the upper
	half-plane $\Hup$. Kalaj and Pavlovi\'c proved that a
	homeomorphism $h$ belongs to $\HQS(\R)$ if and only if $h$ is
	bi-Lipschitz and the Hilbert transform of $h'$ belongs to $L^\infty(\R)$,
	and asked whether $\HQS(\R)$ is a group under composition.
	We give a negative answer to this question. Further, we  show that  the Beltrami coefficients of harmonic quasiconformal mappings of $\Hup$  can be identified with the open unit ball of $H^\infty(\Hup)$.
	Finally, we study the topology of $\HQS(\R)$ and prove that
	it is contractible in the topology inherited from
	the BMO-Teichm\"uller space.
\end{abstract}

\maketitle

\section{Introduction}\label{sec:intro}
Harmonic mappings and quasiconformal mappings are two central classes
of planar mappings. Harmonic mappings arise from the Dirichlet problem
and the theory of minimal energy, while quasiconformal mappings control
geometric distortion and play a basic role in the Teichm\"uller space theory.
The classical Rad\'o--Kneser--Choquet theorem states that the harmonic extension of an orientation-preserving homeomorphism of the unit circle  is a  diffeomorphism of the unit disk (see, for example, \cite{Duren}). 
Imposing quasiconformality on the extension leads to a more restrictive class.  Recall that an increasing homeomorphism $h:\R\to\R$ is called
\emph{quasisymmetric} if there is a constant $M\geq1$ such that
\begin{equation}\label{eq:QS-intro}
	\frac1M\leq\frac{h(x+t)-h(x)}{h(x)-h(x-t)}\leq M
	\qquad(x\in\R,\ t>0).
\end{equation}
The fundamental theorem of Beurling and Ahlfors~\cite{BA} states that
$h$ is quasisymmetric if and only if it is the boundary value of a
quasiconformal mapping of $\Hup$.
Let $\HQC(\Hup)$ be the class of harmonic quasiconformal
homeomorphisms of $\Hup$ fixing $\infty$, and set
$\HQS(\R)=\{F|_{\R}:F\in\HQC(\Hup)\}$.  This raises some natural questions: how to characterize $\HQS(\R)$ and the Beltrami coefficients of the class $\HQC(\Hup)$, as well as their associated algebraic and topological structures. 

Kalaj and Pavlovi\'c~\cite[Theorem~1.2]{KP} proved that an increasing
homeomorphism $h$ belongs to $\HQS(\R)$ if and only if $h$ is
bi-Lipschitz and $\HR(h')\in L^\infty(\R)$. Here $\HR$ denotes the
normalized Hilbert transform on the real line,
\begin{equation}\label{eq:HR}
	\HR f(x)=\frac1\pi\operatorname{p.v.}\!\int_{\R}
	\left(\frac1{x-t}+\frac{t}{1+t^2}\right)f(t)\,dt.
\end{equation} (The case of the unit disk is treated in \cite {Pavlovic}.) It is well known that the Hilbert transform is a bounded linear operator from $L^\infty(\R)$ to $\mathrm{BMO}(\R)$.

In Section \ref{sec:beltrami}, we characterize the Beltrami coefficients of $\HQC(\mathbb{H})$ via the open unit ball of $H^\infty(\mathbb{H})$. Since $\QS(\mathbb{R})$ is a group, Kalaj and Pavlovi\'c naturally asked whether the same holds for $\HQS(\mathbb{R})$. We give a negative answer in Sections \ref{sec:real} and \ref{sec:circle}, and prove in \ref{sec:topology} that the class $\HQS(\mathbb{R})$ is contractible in the topology inherited from the BMO-Teichm\"uller space.
Finally, we  raise some open questions concerning the class \(\HQS(\mathbb{R})\) and chord‑arc curves.

\section{Beltrami coefficients of the class $\HQC(\Hup)$}
\label{sec:beltrami}

 For a quasiconformal mapping $F$, let \(\mu_F=F_{\bar z}/F_z\) denote its Beltrami coefficient. 
Recall that the measurable Riemann mapping theorem
assigns to every $\mu\in L^\infty(\Hup)$ with $\|\mu\|_\infty<1$
a unique quasiconformal homeomorphism
$F_\mu:\Hup\to\Hup$ having Beltrami coefficient $\mu$ almost
everywhere and normalized by
$F_\mu(0)=0$, $F_\mu(1)=1$, and $F_\mu(\infty)=\infty$. By convention, we define
\begin{equation}\label{eq:omega-from-mu}
	\omega_\mu(z)=\overline{\mu(z)}
	\frac{1+\mu(z)}{1+\overline{\mu(z)}}.
\end{equation}

Recall that Kalaj and Pavlovi\'c~\cite[Theorem~3.1]{KP} represented
every  map $F\in \HQC(\Hup)$ in the form
\begin{equation}\label{eq:KP-representation}
	F(x+iy)=u(x,y)+icy,\qquad c>0,
\end{equation}
where $u$ is real-valued and harmonic,
$u_x-iu_y$ is holomorphic and bounded, and $u_x$ has a uniform
positive lower bound. 
The following result is an algebraic reformulation of \eqref{eq:KP-representation}, but  of independent interest.

\begin{theorem}
\label{thm:beltrami}
Let $\mu\in L^\infty(\Hup)$ with $\|\mu\|_\infty<1$,
and let $F_\mu$ be the normalized quasiconformal solution described
above.
Then $F_\mu$ is harmonic if and only if $\omega_\mu$ agrees almost
everywhere with a function in the open unit ball of $H^\infty(\Hup)$.
\end{theorem}
\begin{proof}
Let $F=u+icy\in\HQC(\Hup)$, as in~\eqref{eq:KP-representation},
and define $ q(z)=\frac{1}{c}(u_x(z)-iu_y(z)). $
Then $q$ is bounded and holomorphic, and its image is contained in a
compact subset of the right half-plane. Direct calculation gives
\begin{equation}\label{eq:mu-q}
 F_z=\frac c2(q+1),\qquad
 F_{\bar z}=\frac c2(\overline q-1),\qquad
 \mu_F=\frac{\overline q-1}{q+1}.
\end{equation}
The Cayley transform $ \omega=\frac{q-1}{q+1} $
belongs to $H^\infty(\Hup)$ and satisfies $\|\omega\|_\infty<1$.
Using~\eqref{eq:mu-q} gives
\[
 \mu_F=\overline\omega\frac{1-\omega}{1-\overline\omega},
 \qquad
 \omega=\overline{\mu_F}\frac{1+\mu_F}{1+\overline{\mu_F}}.
\]
This proves necessity.

Conversely, suppose $\omega=\omega_\mu$ belongs to the open unit
ball of $H^\infty(\Hup)$, and write $k=\|\omega\|_\infty<1$. Set $  q_\omega=\frac{1+\omega}{1-\omega}.$
It is bounded and holomorphic, with
\[
 \frac{1-k}{1+k}\leq\operatorname{Re}q_\omega(z)
 \leq\frac{1+k}{1-k}\qquad(z\in\Hup).
\]
Since $\Hup$ is simply connected, the formula
\begin{equation}\label{eq:F-omega}
 F_\omega(z)=\operatorname{Re}\int_i^z q_\omega(\zeta)\,d\zeta
             +i\operatorname{Im}z
\end{equation}
defines a harmonic map. If $u_\omega=\operatorname{Re}\int_i^z
q_\omega(\zeta)\,d\zeta$, then
\[
 (u_\omega)_x=\operatorname{Re}q_\omega>0,
 \qquad (u_\omega)_y=-\operatorname{Im}q_\omega.
\]
For each $y>0$, the map $x\mapsto u_\omega(x,y)$ is an increasing
homeomorphism of $\R$ onto itself, since its derivative has a uniform
positive lower bound. Hence $F_\omega$ is a homeomorphism of
$\Hup$. Moreover, its derivative is bounded and its Jacobian is
$J_{F_\omega}=(u_\omega)_x$ with a uniform positive lower bound;
thus $F_\omega$ is quasiconformal and fixes $\infty$.
Finally,~\eqref{eq:mu-q} gives
$\mu_{F_\omega}=\overline\omega \frac{1-\omega }{1-\overline\omega}=\mu$. Let
$$
\widetilde F_\omega(z)
=\frac{F_\omega(z)-F_\omega(0)}
{F_\omega(1)-F_\omega(0)}.
$$
We obtain a normalized harmonic quasiconformal homeomorphism  with Beltrami coefficient $\mu$. Hence $\widetilde F_\omega=F_\mu$ by uniqueness.
\end{proof}

Note that $\mu$ itself need not belong to $H^\infty(\mathbb{H})$; in fact, requiring $\mu\in H^\infty(\mathbb{H})$ forces the map to be affine.

\begin{corollary}
\label{cor:analytic-mu}
Let $a\in H^\infty(\Hup)$ with $\|a\|_\infty<1$, and let $F_a$  be the normalized quasiconformal solution described
above. 
Then $F_a$ is harmonic if and only if $a$ is constant.
\end{corollary}
\begin{proof}
If $F_a$ is harmonic, Theorem~\ref{thm:beltrami} implies that
$\omega_a=\overline a(1+a)/(1+\overline a)$ is holomorphic. Since $a$
is holomorphic, $ \bar\partial\omega_a
=\frac{1+a}{(1+\overline a)^2}\,\overline{a'}. $
The factor preceding $\overline{a'}$ does not vanish, so this occurs
only when $a'\equiv0$. Conversely, if $a$ is constant, the affine map $F_a(z)=\frac{z+a\overline z}{1+a} $. 
is a harmonic quasiconformal homeomorphism of $\Hup$ fixing
$0,1,\infty$ and having Beltrami coefficient $a$. Uniqueness of the
normalized solution proves the assertion.
\end{proof}
\begin{remark}
	In fact, this raises another interesting question: what special properties does the Beltrami solution \(F_\mu\) possess when \(\mu\in H^\infty(\mathbb H)\)? We prove that  \(F_\mu\) is bi-Lipschitz, real-analytic, and fixes the boundary pointwise; see \cite{LiuShen}.
\end{remark}

For a bounded interval $I\subset\R$, let
$Q_I=I\times(0,|I|)\subset\Hup$. Recall that a positive measure
$\lambda$ on $\Hup$ is a \emph{Carleson measure} if $\sup_{I\subset\R}\frac{\lambda(Q_I)}{|I|}<\infty. $
We also compare these coefficients with the Carleson-measure
condition appearing in the BMO-Teichm\"uller space theory. 

\begin{proposition}
\label{prop:carleson-rigidity}
Let $F\in\HQC(\Hup)$. Then $  \frac{|\mu_F(z)|^2}{\operatorname{Im}z}\,dx\,dy$
is a Carleson measure on $\Hup$ if and only if $F$ is conformal.
More precisely, if $F$ is not conformal, then for every 
bounded interval $I\subset\R$,
\begin{equation}\label{eq:infinite-carleson}
 \int_0^{|I|}\int_I
      \frac{|\mu_F(x+iy)|^2}{y}\,dx\,dy=+\infty.
\end{equation}
\end{proposition}
\begin{proof}
Take $\omega\in H^\infty(\Hup)$ corresponding to $\mu_F$ in
Theorem~\ref{thm:beltrami}. Then $|\mu_F|=|\omega|$. If $F$ is
nonconformal, $\omega\not\equiv0$. By Fatou's theorem, $\omega$ has
nontangential boundary values $\omega^*(x)$ almost everywhere; by the
boundary uniqueness theorem for bounded analytic functions,
$\omega^*(x)\neq0$ almost everywhere on $\R$; see~\cite{Garnett}.
Consequently, for almost every $x\in I$,
\[
 \int_0^{|I|}\frac{|\omega(x+iy)|^2}{y}\,dy=+\infty,
\]
since $\omega(x+iy)\to\omega^*(x)\neq0$ as $y\downarrow0$.
This proves~\eqref{eq:infinite-carleson}. If $F$ is
conformal, then $\mu_F\equiv0$, so the measure vanishes.
\end{proof}

\section{The real-line counterexample}\label{sec:real}
Recall  the boundary characterization of Kalaj and Pavlovi\'c: for an increasing homeomorphism $h$,
\begin{equation}\label{eq:KP}
 h\in\HQS(\R)
 \quad\Longleftrightarrow\quad
 h\text{ is bi-Lipschitz and }\HR(h')\in L^\infty(\R).
\end{equation}
Our counterexample is constructed directly at the level of boundary
values.

\begin{theorem}\label{thm:line}
There exists $\varepsilon_0>0$ and an explicit family $g_\varepsilon\in\HQS(\R)$, $0<\varepsilon<\varepsilon_0$, with
\[
 g_\varepsilon(0)=0,\qquad g_\varepsilon(1)=1,
\]
such that
\[
 g_\varepsilon\circ g_\varepsilon\notin\HQS(\R),
 \qquad g_\varepsilon^{-1}\notin\HQS(\R).
\]
Moreover,
\begin{equation}\label{eq:smallline}
 \|g_\varepsilon'-1\|_\infty+
 \|\HR(g_\varepsilon'-1)\|_\infty=O(\varepsilon)
 \qquad(\varepsilon\downarrow0).
\end{equation}
\end{theorem}

The obstruction is the following compatibility of logarithmic means. We use the classical fact that, if $f,\HR f\in L^\infty(\R)$ and $f$ is real-valued, then $f$ is the boundary real part of a function in $H^\infty(\Hup)$; see~\cite{Garnett}.

\begin{lemma}\label{lem:mean}
Let $a>1$ and let $f\in L^\infty(\R)$ be real-valued, with $f(ax)=f(x)$ almost everywhere. If $\HR f\in L^\infty(\R)$, then
\begin{equation}\label{eq:means}
 \frac1{\log a}\int_0^{\log a}f(e^s)\,ds
 =\frac1{\log a}\int_0^{\log a}f(-e^s)\,ds.
\end{equation}
\end{lemma}
\begin{proof}
Choose $F(z)\in H^\infty(\Hup)$ such that  Re$F|_{\R}=f$. Since $F(az)-F(z)$ has boundary real part zero, it is a purely imaginary constant $i\beta$. Iteration  gives  $F(a^nz)=F(z)+ni\beta$, and hence $\beta=0$. Thus $F(e^w)$ is bounded and analytic on $0<\operatorname{Im}w<\pi$, and is periodic with period $L=\log a$. Writing $F(e^w)=A(x,y)+iB(x,y)$, the Cauchy--Riemann equations give
\[
 \frac{d}{dy}\frac1L\int_0^L A(x, y)\,dx
 =-\frac1L\int_0^L B_x(x, y)\,dx=0.
\]
Taking boundary limits at $y=0$ and $y=\pi$ and using bounded convergence proves~\eqref{eq:means}.
\end{proof}

\begin{proof}[Proof of Theorem~\ref{thm:line}]
For the principal branch of $\operatorname{Log}$ on $\Hup$, put
\[
 \Phi(z)=z^i=e^{i\operatorname{Log}z},\qquad \rho=e^{-\pi}.
\]
Since $|\Phi(z)|=e^{-\arg z}\leq1$, its boundary real part is
\begin{equation}\label{eq:u}
 u(x)=\begin{cases}
 \cos(\log x),&x>0,\\
 \rho\cos(\log|x|),&x<0,
 \end{cases}
\end{equation}
and $u,\HR u\in L^\infty(\R)$. For $0<\varepsilon<1$, define
\begin{equation}\label{eq:g}
 U(x)=\int_0^x u(t)\,dt,\qquad
 C_\varepsilon=\frac1{1+\varepsilon/2},\qquad
 g_\varepsilon(x)=C_\varepsilon\bigl(x+\varepsilon U(x)\bigr).
\end{equation}
Since $U(0)=0$ and $U(1)=1/2$, we have $g_\varepsilon(0)=0$ and $g_\varepsilon(1)=1$. Moreover,
\[
 C_\varepsilon(1-\varepsilon)\leq g_\varepsilon'
 \leq C_\varepsilon(1+\varepsilon)\quad\text{a.e.},\qquad
 g_\varepsilon'-1=C_\varepsilon\varepsilon\left(u-\tfrac12\right).
\]
As $u,\HR u\in L^\infty(\R)$ and $\HR 1=0$, criterion~\eqref{eq:KP} gives $g_\varepsilon\in\HQS(\R)$, together with~\eqref{eq:smallline}. Also,
\[
 g_\varepsilon(e^{2\pi}x)=e^{2\pi}g_\varepsilon(x).
\]
Thus $d_\varepsilon=(g_\varepsilon\circ g_\varepsilon)'$ is multiplicatively periodic with factor $e^{2\pi}$.

Set
\begin{equation}\label{eq:B}
 B_p(s)=1+\frac p2(\cos s+\sin s),\qquad
\end{equation}
The identity $\int_0^r\cos(\log t)\,dt=\frac r2(\cos\log r+\sin\log r)$ yields
\[
 g_\varepsilon(e^s)=C_\varepsilon e^s B_\varepsilon(s),\qquad
 g_\varepsilon(-e^s)=-C_\varepsilon e^s B_{\rho\varepsilon}(s).
\]
Thus $d_\varepsilon(e^s)=D_{\varepsilon,\varepsilon}(s)$ and
$d_\varepsilon(-e^s)=D_{\rho\varepsilon,\varepsilon}(s)$, where
\begin{equation}\label{eq:D}
 D_{p,\varepsilon}(s)=C_\varepsilon^2
 \bigl[1+p\cos(s+\log C_\varepsilon+\log B_p(s))\bigr](1+p\cos s).
\end{equation}
Uniformly in $s$ and for $\log C_\varepsilon$ in a fixed neighborhood of $0$, as $p\to0$,
\begin{align*}
 C_\varepsilon^{-2}D_{p,\varepsilon}(s)
 &=1+p\bigl(\cos(s+\log C_\varepsilon)+\cos s\bigr)\\
 &\quad+p^2\left[\cos(s+\log C_\varepsilon)\cos s
 -\tfrac12(\cos s+\sin s)\sin(s+\log C_\varepsilon)\right]+O(p^3).
\end{align*}
The period mean is consequently
\begin{equation}\label{eq:M}
 M_\varepsilon(p):=\frac1{2\pi}\int_0^{2\pi}D_{p,\varepsilon}(s)\,ds
 =C_\varepsilon^2\left[1+\frac{p^2}{4}
 (\cos\log C_\varepsilon-\sin\log C_\varepsilon)+O(p^3)\right].
\end{equation}
Since $C_\varepsilon=1+O(\varepsilon)$ and $\log C_\varepsilon=O(\varepsilon)$, we obtain
\begin{equation}\label{eq:gap}
 M_\varepsilon(\varepsilon)-M_\varepsilon(\rho\varepsilon)
 =\frac{1-\rho^2}{4}\varepsilon^2+O(\varepsilon^3)>0
\end{equation}
for sufficiently small $\varepsilon>0$. Lemma~\ref{lem:mean} implies $\HR d_\varepsilon\notin L^\infty(\R)$, proving the assertion for self-composition.

For the inverse, on either half-axis put $t=s+\log C_\varepsilon+\log B_p(s)$, with $p=\varepsilon$ or $p=\rho\varepsilon$, respectively. Then
\[
 \frac{dt}{ds}=\frac{1+p\cos s}{B_p(s)},\qquad
 (g_\varepsilon^{-1})'(\pm e^t)=\frac1{C_\varepsilon(1+p\cos s)}.
\]
The latter function is $2\pi$-periodic in $t$. Its period mean is
\begin{equation}\label{eq:N}
 N_\varepsilon(p)=\frac1{2\pi C_\varepsilon}\int_0^{2\pi}\frac{ds}{B_p(s)}
 =\frac1{C_\varepsilon\sqrt{1-p^2/2}}.
\end{equation}
Since $N_\varepsilon(\varepsilon)\neq N_\varepsilon(\rho\varepsilon)$, another application of Lemma~\ref{lem:mean} and~\eqref{eq:KP} gives $g_\varepsilon^{-1}\notin\HQS(\R)$.
\end{proof}

\begin{remark}\label{rem:extensions}
The inverse of a boundary map must be distinguished from the inverse of a particular harmonic extension. In fact, harmonic quasiconformal extension from $\R$ is not unique: for any $a\in\R$ and $c>0$, $ x+ay+icy $
is a harmonic quasiconformal homeomorphism of $\Hup$ with identity boundary values. Thus one cannot infer the failure of a harmonic extension for $h^{-1}$ merely by observing that the interior inverse of a chosen extension of $h$ is not harmonic. 
\end{remark}

\section{The circle counterexample}\label{sec:circle}
We now turn to the unit disk $\D$ and its boundary $\T=\partial\D$.    A direct M\"obius conjugation of its harmonic extension would not suffice. If $F$ is harmonic and \(\varphi\) is conformal,  \(F\circ \varphi\) need not be harmonic, whereas  \(\varphi\circ F\) preserves harmonicity.   For this reason, a technical adaptation of the real-line argument is required for the circle setting.

Let $\HQS(\T)$ consist of the orientation-preserving circle homeomorphisms whose Poisson extensions are quasiconformal. Write $h(e^{it})=e^{i\varphi(t)}$, where $\varphi$ is increasing and $\varphi(t+2\pi)=\varphi(t)+2\pi$. Pavlovi\'c's theorem~\cite[Theorem~1.1]{Pavlovic} says that
\begin{equation}\label{eq:Pav}
 h\in\HQS(\T)
 \quad\Longleftrightarrow\quad
 \varphi\text{ is bi-Lipschitz and }\HT(\varphi')\in L^\infty(\T),
\end{equation}
where
\[\HT v(\theta)=\frac1{2\pi}\operatorname{p.v.}\!\int_{-\pi}^{\pi}
 v(t)\cot\frac{\theta-t}{2}\,dt.
\]

\begin{lemma}\label{lem:local}
Suppose $v\in L^\infty(\T)$ is real-valued and, for some $0<r_0<\pi$,
\[
 v(t)=V_+(\log t),\qquad v(-t)=V_-(\log t)
 \quad(0<t<r_0),
\]
where $V_\pm$ are bounded periodic functions of a common period. If their period means $m_\pm$ are different, then $\HT v\notin L^\infty(\T)$.
\end{lemma}
\begin{proof}
Let $Q_r(t)=2r\sin t/(1-2r\cos t+r^2)$ be the conjugate Poisson kernel and put $a=1-r$. Uniformly as $r\uparrow1$,
\[
 \int_0^{r_0}\left|Q_r(t)-\frac2t\mathbf1_{(a,r_0)}(t)\right|dt=O(1).
\]
Indeed, the integral on $(0,a)$ is bounded, while on $(a,r_0)$ the difference is $O(t+a^2/t^3)$. The remaining part of the circle is separated from the singularity. Consequently,
\begin{align*}
 \frac1{2\pi}\int_{-\pi}^{\pi}Q_r(-t)v(t)\,dt
 &=\frac1\pi\int_a^{r_0}\frac{v(-t)-v(t)}{t}\,dt+O(1)\\
 &=\frac1\pi\int_{\log a}^{\log r_0}(V_-(s)-V_+(s))\,ds+O(1)\\
 &=\frac{m_--m_+}{\pi}\log\frac{r_0}{a}+O(1).
\end{align*}
If $\HT v$ were bounded,   the conjugate Poisson integral of $v$  would be uniformly bounded. This is a contradiction as $a = 1 - r \to 0^+$.

\end{proof}

\begin{theorem}\label{thm:circle}
There are $h_\varepsilon=e^{i\phi_{\varepsilon}}\in\HQS(\T)$, equal to the identity outside a fixed small arc about $1$, such that
\[
 h_\varepsilon\circ h_\varepsilon\notin\HQS(\T),
 \qquad h_\varepsilon^{-1}\notin\HQS(\T),
\]
\begin{equation}\label{eq:smallcircle}
 \|\varphi_\varepsilon'-1\|_\infty+
 \|\HT(\varphi_\varepsilon'-1)\|_\infty=O(\varepsilon).
\end{equation}
\end{theorem}
\begin{proof}
Choose a real $\chi\in C_c^\infty(-\pi/4,\pi/4)$ equal to $1$ near $0$. With $U$ and $C_\varepsilon$ as in~\eqref{eq:g}, define on $[-\pi,\pi]$
\begin{equation}\label{eq:cutoff}
 \varphi_\varepsilon(t)=t+\varepsilon C_\varepsilon\chi(t)
 \left(U(t)-\frac t2\right),
\end{equation}
and extend by $\varphi_\varepsilon(t+2\pi)=\varphi_\varepsilon(t)+2\pi$. Then
\[
 \varphi_\varepsilon'=1+\varepsilon C_\varepsilon b,
 \qquad b=\chi\left(u-\tfrac12\right)
 +\chi'\left(U-\tfrac t2\right).
\]

We claim that $b, \HT b\in L^\infty(\T)$. Clearly $b\in L^\infty(\T)$, and the terms $\chi'(U-t/2)$ and $-\frac12 \varepsilon C_\varepsilon\chi$ are smooth on $\T$, so their circle Hilbert transforms are bounded. Moreover, $\HR(\chi u)\in L^\infty(\R)$ because $\chi$ is compactly supported and smooth, and the commutator
\[
\HR(\chi u)(x)-\chi(x)\HR u(x)=\frac1\pi\int_\R\bigl(\chi(t)-\chi(x)\bigr)\left(\frac1{x-t}+\frac{t}{1+t^2}\right)u(t)\,dt
\]
is bounded on $\R$. Since the kernel difference $\frac12\cot\frac{x-t}{2}-\frac1{x-t}$ is bounded for $x\in[-\pi,\pi]$ and $t\in\operatorname{supp}\chi$, it follows that $\HT(\chi u)\in L^\infty(\T)$, proving the claim.

For small $\varepsilon>0$, $\varphi_\varepsilon$ is bi-Lipschitz. Thus~\eqref{eq:Pav} gives
$h_\varepsilon(e^{it})=e^{i\varphi_\varepsilon(t)}\in\HQS(\T)$, and~\eqref{eq:smallcircle} holds. The normalization follows from~\eqref{eq:cutoff}.

Near $0$, $\varphi_\varepsilon=g_\varepsilon$, and, after shrinking the neighborhood,
\[
 \varphi_\varepsilon\circ\varphi_\varepsilon
 =g_\varepsilon\circ g_\varepsilon,
 \qquad \varphi_\varepsilon^{-1}=g_\varepsilon^{-1}.
\]
Their derivatives therefore have the same unequal local logarithmic means as in~\eqref{eq:gap} and~\eqref{eq:N}, respectively. The proof is completed by Lemma~\ref{lem:local} and~\eqref{eq:Pav}.
\end{proof}

\section{Contractibility in the BMO-Teichm\"uller space}\label{sec:topology}

For a locally integrable real-valued function $u$ on $\R$, write
\begin{equation}\label{eq:BMO}
 \|u\|_{\BMO(\R)}
 =\sup_{I\subset\R}\frac1{|I|}\int_I|u(x)-u_I|\,dx,
 \qquad u_I=\frac1{|I|}\int_Iu(x)\,dx,
\end{equation}
where $I$ ranges over all bounded intervals. The space $\BMO(\R)$ consists
of functions for which this quantity is finite; as usual, its elements are
understood modulo additive constants.

An increasing homeomorphism $h:\R\to\R$ is called \emph{strongly
	quasisymmetric} in the sense of Semmes~\cite{Semmes}, written
$h\in\SQS(\R)$, if there exist $C_1,C_2>0$ such that
\begin{equation}\label{eq:SQS}
	\frac{|h(E)|}{|h(I)|}
	\le C_1\left(\frac{|E|}{|I|}\right)^{C_2}
\end{equation}
for every interval $I\subset\R$ and every measurable set $E\subset I$.
Equivalently, $h$ is locally absolutely continuous and $h'$ is an
$A_\infty$-weight  introduced by Muckenhoupt (see \cite{Garnett,Muck}). In particular, $\log h'\in\BMO(\R)$.
 The class $\SQS(\R)$ forms a subgroup of $\QS(\R)$. Moreover,
 $h\in\SQS(\R)$ if and only if it admits a quasiconformal extension
 to $\Hup$ with Beltrami coefficient $\mu$ such that $\frac{|\mu(z)|^2}{\operatorname{Im}z}\,dx\,dy $
 is a Carleson measure; see \cite{FKP,Semmes}. Every bi-Lipschitz
increasing homeomorphism is strongly quasisymmetric, so $ \HQS(\R)\subset\SQS(\R).$
For $h\in\HQS(\R)$, however, the quasiconformal extension
satisfying the Carleson condition cannot be harmonic, by
Proposition~\ref{prop:carleson-rigidity}.

Let $\operatorname{Aff}^+(\R)$ be the group of increasing affine maps
$x\mapsto ax+b$ with $a>0$. The real-line model of the
\emph{BMO-Teichm\"uller space} is $  T_b(\R)=\SQS(\R)/\operatorname{Aff}^+(\R).$
Every equivalence class has a unique representative fixing $0$ and $1$.
The logarithmic-derivative coordinate $[h]\longmapsto[\log h']\in\BMO(\R) $
realizes the BMO-Teichm\"uller topology.  This class  first appeared in \cite{AZ}, where it is also shown to be path‑connected.
Further detailed developments may be found in  \cite{FHS,LiuShenVMO, Shen,SW}

In this section and the next, all maps in $\HQS(\R)$ are
normalized by $h(0)=0$ and $h(1)=1$. We use the metric
\begin{equation}\label{eq:dB}
	d(h,g)=\|\log h'-\log g'\|_{\BMO(\R)}.
\end{equation}

It is easy to see $\HQS(\R)$ is convex. For finitely many $h_j\in\HQS(\R)$, choose
 $0<m\le h_j'\le M$ a.e. If $\lambda_j\ge0$ and
 $\sum_j\lambda_j=1$, then
 \[
 m\le\big(\sum_j\lambda_jh_j\big)'\le M,
 \qquad
 \HR\!\big(\sum_j\lambda_jh_j\big)'
 =\sum_j\lambda_j\HR(h_j')\in L^\infty(\R).
 \]  
 When $\HQS(\R)$  is placed in the natural BMO–Teichm\"uller setting, its topology properties and its relation to the well-known chord-arc curves are our main concerns.

 \begin{proposition}\label{prop:not-open-closed}
 	The normalized class $\HQS(\R)$ is neither open nor closed in
 	$T_b(\R)$.
 \end{proposition}
 \begin{proof}
 	For $\varepsilon>0$, set
 	\[
 	h_\varepsilon(x)=\operatorname{sgn}(x)|x|^{1+\varepsilon},
 	\qquad
 	\log h_\varepsilon'=\log(1+\varepsilon)+\varepsilon\log|x|.
 	\]
 	Since $|x|^\varepsilon\in A_\infty$ but $h_\varepsilon$ is not
 	bi-Lipschitz, we have
 	$h_\varepsilon\in\SQS(\R)\setminus\HQS(\R)$, while
 	\[
 	d(h_\varepsilon,\id)
 	=\varepsilon\|\log|x|\|_{\BMO(\R)}\longrightarrow0.
 	\]
 	Thus $\HQS(\R)$ is not open. More generally, normalizing $g_\epsilon(x)=\int_0^x h'(t)|t|^\epsilon\,dt$ for any $h\in \HQS(\R)$ shows that $\HQS(\R)$ has empty interior in $T_b(\R)$.
 	
 	For nonclosedness, choose $0\le\chi\in C_c^\infty((-1/2,1/2))$
 	with $\chi=1$ on $[-1/4,1/4]$, and put
 	\[
 	v(x)=\begin{cases}
 		\displaystyle\frac{\chi(x)}{\log(e/x)},&x>0,\\[2mm]
 		0,&x\le0.
 	\end{cases}
 	\]
 	Choose $0\le v_n\in C_c^\infty(\R)$ with $v_n\to v$
 	uniformly and define
 	\[
 	h_n(x)=\frac{\int_0^x(1+v_n(t))\,dt}
 	{\int_0^1(1+v_n(t))\,dt},
 	\qquad
 	h(x)=\frac{\int_0^x(1+v(t))\,dt}
 	{\int_0^1(1+v(t))\,dt}.
 	\]
 	Since $\HR v_n\in L^\infty(\R)$ and $1+v_n\ge1$, we have
 	$h_n\in\HQS(\R)$ and $h\in\SQS(\R)$. Moreover,
 	\[
 	d(h_n,h)\le
 	2\|\log h_n'-\log h'\|_\infty\longrightarrow0.
 	\]
 	On the other hand, for $a\downarrow0$,
 	\[
 	-\HR v(-a)
 	=\frac1\pi\int_0^{1/4}\frac{dt}{(a+t)\log(e/t)}+O(1)
 	\ge\frac1{2\pi}\int_a^{1/4}\frac{dt}{t\log(e/t)}+O(1)
 	\longrightarrow\infty.
 	\]
 	Hence $\HR(h')\notin L^\infty(\R)$, so
 	$h\notin\HQS(\R)$.
 \end{proof}

\begin{theorem}\label{thm:contract}
	The normalized class $\HQS(\R)$ is contractible in the topology inherited
	from $T_b(\R)$. More precisely, it admits a strong deformation
	retraction onto the identity.
\end{theorem}

To prove this result, we need the following lemma.

\begin{lemma}\label{lem:logconvex}
There exist constants $\delta_0>0$ and $C\geq 1$ such that for any 
$u_1,\ldots,u_N\in\BMO(\R)$ with $\delta:=\max_{1\le j\le N}\|u_j\|_{\BMO(\R)}\leq\delta_0$ 
and any $\lambda_1,\ldots,\lambda_N\geq 0$ with $\sum_{j=1}^N\lambda_j=1$,
\begin{equation}\label{eq:logconvex}
	\|\log\sum_{j=1}^N\lambda_j e^{u_j}\|_{\BMO(\R)}\leq C\delta,
\end{equation}
where $\delta_0$ and $C$ are independent of $N$.
\end{lemma}

\begin{proof}
	The case $\delta=0$ is immediate, so we assume $0<\delta\leq\delta_0$. 
	Fix a bounded interval $I\subset\R$. 
	By  John--Nirenberg inequality~\cite[Chapter~VI]{Garnett}, there exist universal constants $C_1,C_2>0$ such that for each $1\leq j\leq N$,
	\[
	\frac{1}{|I|}\bigl|\{x\in I:|u_j(x)-(u_j)_I|>t\}\bigr|
	\leq C_1e^{-\frac{C_2}{\delta}t}\qquad\forall t>0.
	\]
	Taking $\delta_0:=C_2/2$, so that $C_2/\delta-1\geq C_2/(2\delta)>0$, we have 
	\begin{equation}\label{eq:expmoment}
		\begin{aligned}
			\frac{1}{|I|}\int_I \bigl(e^{|u_j(x)-(u_j)_I|}-1\bigr)\,dx
			&=\frac{1}{|I|}\int_0^\infty e^t\,\bigl|\{x\in I:|u_j(x)-(u_j)_I|>t\}\bigr|\,dt\\
			&\leq C_1\int_0^\infty e^{-(C_2/\delta-1)t}\,dt\\
			&\leq \frac{2C_1}{C_2}\delta.
		\end{aligned}
	\end{equation}
	In particular, all exponentials in~\eqref{eq:logconvex} are locally integrable. 
	
Define
	\[ U:=\log\sum_{j=1}^N\lambda_j e^{u_j}, \quad
	c_I:=\log\sum_{j=1}^N\lambda_j e^{(u_j)_I},
	\qquad
	\alpha_j:=\frac{\lambda_j e^{(u_j)_I}}{\sum_{k=1}^N\lambda_k e^{(u_k)_I}}.
	\]
	Then $\alpha_j\geq 0$, $\sum_{j=1}^N\alpha_j=1$, and $v:=U-c_I=\log\sum_{j=1}^N\alpha_j e^{u_j-(u_j)_I}$. 
	Applying Jensen's inequality for the lower bound and the inequality $\log t\leq t-1$ for the upper bound, we get
	\[
	\sum_{j=1}^N\alpha_j\bigl(u_j-(u_j)_I\bigr)
	\leq v
	\leq \sum_{j=1}^N\alpha_j e^{u_j-(u_j)_I}-1
	=\sum_{j=1}^N\alpha_j\bigl(e^{u_j-(u_j)_I}-1\bigr).
	\]
	Since $ |e^t-1|\leq e^{|t|}-1$ for all $t\in\R$, this implies
	\[
	|v|\leq\sum_{j=1}^N\alpha_j\bigl(e^{|u_j-(u_j)_I|}-1\bigr).
	\]
	Then
	\[
	\frac{1}{|I|}\int_I |U(x)-(U)_I|\,dx
	\leq \frac{2}{|I|}\int_I |v(x)|\,dx
	\leq \frac{4C_1}{C_2}\delta.
	\]
	Taking the supremum over all  $I\subset\R$ proves~\eqref{eq:logconvex} with $C:=\max\{1,\frac{4C_1}{C_2}\}$.
\end{proof}

\begin{proof}[Proof of Theorem~\ref{thm:contract}]	
	We first note that finite convex combinations in $\HQS(\R)$ depend continuously on their coefficients. Indeed, given fixed $h_1,\ldots,h_N\in\HQS(\R)$, there exist constants $0<m\leq M<\infty$ such that $m\leq h_j'\leq M$ almost everywhere for each $1\leq j\leq N$. For any convex combinations $h_\lambda=\sum_{j=1}^N\lambda_jh_j$ and $h_\mu=\sum_{j=1}^N\mu_jh_j$, the mean value theorem yields
	\begin{equation}\label{eq:simplex}
		d(h_\lambda,h_\mu)
		\leq 2\|\log h_\lambda'-\log h_\mu'\|_\infty
		\leq \frac{2M}{m}\sum_{j=1}^N|\lambda_j-\mu_j|.
	\end{equation}
	
	Moreover, we show that small balls have uniformly small convex hulls. Fix $h\in\HQS(\R)$ and $0<r\leq\delta_0$. Consider any finite convex combination $h_\lambda = \sum_{j=1}^N \lambda_j h_j$ of elements $h_1,\ldots,h_N \in B_d(h,r)$. To bound the distance $d(h_\lambda, h)$, we set $u_j = \log h_j' - \log h'$. Since $h_j \in B_d(h,r)$, we have $\|u_j\|_{\BMO(\R)} = d(h_j, h) < r \leq \delta_0$. Applying Lemma~\ref{lem:logconvex} to these $u_j$, we obtain
	\begin{align*}
		d(h_\lambda, h) &= \| \log \Big( \sum_{j=1}^N \lambda_j h_j' \Big) - \log h' \|_{\BMO(\R)} \\
		&= \| \log \sum_{j=1}^N \lambda_j e^{u_j}\|_{\BMO(\R)} \\
		&\leq C \max_{1\le j\le N} \|u_j\|_{\BMO(\R)} < Cr.
	\end{align*}
	Thus, $h_\lambda \in B_d(h,Cr)$. This proves that for every $h\in\HQS(\R)$ and $0<r\leq\delta_0$,
	\begin{equation}\label{eq:convball}
		\conv B_d(h,r)\subset B_d(h,Cr),
	\end{equation}
	where $C\geq1$ is the universal constant from Lemma~\ref{lem:logconvex}.
	
	Now choose a decreasing sequence $r_n\in(0,\delta_0]$ such that $Cr_n<2^{-n}$ for all $n\geq1$. For each $n$, consider the open cover $\mathcal U_n$ of $\HQS(\R)$ given by
	\[
	\mathcal U_n = \{B_d(\id,r_n)\} \cup \bigl\{ B_d\bigl(g,\min\{r_n,d(g,\id)/2\}\bigr) : g\in\HQS(\R)\setminus\{\id\} \bigr\}.
	\]
	Since metric spaces are paracompact, there exists a locally finite partition of unity $\{\lambda_{n,j}\}_j$ subordinate to $\mathcal U_n$. For each index $j$, choose a  ball $B_d(p_{n,j},\rho_{n,j})\in\mathcal U_n$ containing $\operatorname{supp}\lambda_{n,j}$. 
	Using these centers $p_{n,j}$ and weight functions $\lambda_{n,j}$, we define the $n$-th approximation map $A_n : \HQS(\R) \to \HQS(\R)$ by
	\[
	A_n(h)=\sum_j\lambda_{n,j}(h)p_{n,j}.
	\]
	By local finiteness and~\eqref{eq:simplex}, each $A_n:\HQS(\R)\to\HQS(\R)$ is continuous. By construction of $\mathcal U_n$, no ball centered at $g\neq \id$ contains $\id$, so $\lambda_{n,j}(\id)>0$ implies $p_{n,j}=\id$; hence $A_n(\id)=\id$. Furthermore, whenever $\lambda_{n,j}(h)>0$, we have $d(h,p_{n,j})<\rho_{n,j}\leq r_n$. It follows from~\eqref{eq:convball} that $A_n(h)\in B_d(h,Cr_n)$ for all $h\in\HQS(\R)$.
	
	Finally, set $t_n=2^{-n}$ ($n\geq1$) and define $\mathcal D:\HQS(\R)\times[0,1]\to\HQS(\R)$ by $\mathcal D(h,0)=h$,
	\[
	\mathcal D(h,t)=\frac{t_n-t}{t_n-t_{n+1}}A_{n+1}(h)+\frac{t-t_{n+1}}{t_n-t_{n+1}}A_n(h) \quad \text{for } t\in[t_{n+1},t_n],
	\]
	and $\mathcal D(h,t)=(2-2t)A_1(h)+(2t-1)\id$ for $t\in[1/2,1]$. Clearly, $\mathcal D(h,1)=\id$ and $\mathcal D(\id,t)=\id$ for all $t\in[0,1]$. On $\HQS(\R)\times(0,1]$, local finiteness of the partitions of unity and~\eqref{eq:simplex} ensure that $\mathcal D$ is jointly continuous. To verify continuity at $t=0$, note that for any $t\in[t_{n+1},t_n]$ and $h\in\HQS(\R)$, the point $\mathcal D(h,t)$ is a finite convex combination of centers $p_{n,j}$ and $p_{n+1,k}$ lying in $B_d(h,r_n)$ (as $r_{n+1}\leq r_n$). Thus~\eqref{eq:convball} gives
	\[
	d(\mathcal D(h,t),h) < Cr_n < 2^{-n} \qquad \forall\, h\in\HQS(\R),\ t\in[t_{n+1},t_n].
	\]
	Consequently, $\mathcal D(h,t)\to h$ as $t\to0^+$ uniformly on $\HQS(\R)$, which implies joint continuity of $\mathcal D$ at $t=0$. This completes the proof.
\end{proof}

\begin{remark}
	The proof relies on  classical topological tools. First, we  construct approximating maps by partition-of-unity convex combinations, adapting the idea from the proof of Dugundji's extension theorem \cite{Dug}.  Second, we  assemble these approximations into a strong deformation retraction via a telescoping homotopy,  a standard technique for  contractibility.
	\end{remark}

\section{Further questions}\label{sec:questions}
Let $T_C(\R)\subset T_b(\R)$ denote the conformal welding class of chord-arc curves, which is an open neighborhood of the identity in $T_b(\R)$ (see \cite{AZ}).  Since the identity is not an interior point of
$\HQS(\R)$ by Proposition~\ref{prop:not-open-closed}, $T_C(\R)\not\subset\HQS(\R). $
However, is it true that $\operatorname{HQS}(\mathbb R)\subset T_C(\mathbb R)$? Furthermore, how can one geometrically characterize quasicircles whose conformal welding belongs to $\operatorname{HQS}(\mathbb R)$?
\section*{Acknowledgments}
Google's Gemini 4 Argon independently constructed the counterexamples  in Sections \ref{sec:real} and \ref{sec:circle}.  The authors assume full responsibility for the paper.

\end{document}